\documentclass[12pt]{article}
\usepackage{amsfonts,amsmath}

\def \beq {\begin{eqnarray}}
\def \eeq {\end{eqnarray}}
\def \beqn {\begin{eqnarray*}}
\def \eeqn {\end{eqnarray*}}

\newcommand{\halmos}{\rule{1ex}{1.4ex}}

\newcounter{for}[section]

\numberwithin{equation}{section}

\newtheorem{itlemma}{Lemma}[section]
\newtheorem{itproposition}[itlemma]{Proposition}
\newtheorem{theorem}[itlemma]{Theorem}
\newtheorem{itcorollary}[itlemma]{Corollary}
\newtheorem{itremark}[itlemma]{Remark}
\newtheorem{itremarks}[itlemma]{Remarks}
\newtheorem{itdefinition}[itlemma]{Definition}
\newtheorem{itexample}[itlemma]{Example}

\newenvironment{fact}{\begin{itfact}\rm}{\end{itfact}}
\newenvironment{claim}{\begin{itclaim}\rm}{\end{itclaim}}
\newenvironment{lemma}{\begin{itlemma}}{\end{itlemma}}
\newenvironment{remark}{\begin{itremark}\rm}{\end{itremark}}
\newenvironment{remarks}{\begin{itremarks} \rm}{\end{itremarks}}
\newenvironment{corollary}{\begin{itcorollary}}{\end{itcorollary}}
\newenvironment{proposition}{\begin{itproposition}}{\end{itproposition}}
\newenvironment{definition}{\begin{itdefinition}\rm}{\end{itdefinition}}
\newenvironment{example}{\begin{itexample}\rm}{\end{itexample}}
\newenvironment{proof}{\noindent {\em Proof}.\ \
}{\hspace*{\fill}$\halmos$\medskip}
\newcommand{\be}[1]{\addtocounter{for}{1} \begin{equation}\label{#1}}
\newcommand{\ee}{\end{equation}}
\newcommand{\bl}[1]{\begin{lemma}\label{#1}}
\newcommand{\br}[1]{\begin{remark}\label{#1}}
\newcommand{\brs}[1]{\begin{remarks}\label{#1}}
\newcommand{\bt}[1]{\begin{theorem}\label{#1}}
\newcommand{\bd}[1]{\begin{definition}\label{#1}}
\newcommand{\bp}[1]{\begin{proposition}\label{#1}}
\newcommand{\bc}[1]{\begin{corollary}\label{#1}}
\newcommand{\bfact}[1]{\begin{fact}\label{#1}}
\newcommand{\bex}[1]{\begin{example}\label{#1}}
\newcommand{\ec}{\end{corollary}}
\newcommand{\efact}{\end{fact}}
\newcommand{\eex}{\end{example}}
\newcommand{\el}{\end{lemma}}
\newcommand{\er}{\end{remark}}
\newcommand{\ers}{\end{remarks}}
\newcommand{\et}{\end{theorem}}
\newcommand{\ed}{\end{definition}}
\newcommand{\ep}{\end{proposition}}
\newcommand{\epr}{\end{proof}}
\newcommand{\bpr}{\begin{proof}}
\newcommand{\bcl}[1]{\begin{claim}\label{#1}}
\newcommand{\ecl}{\end{claim}}

\newcommand{\ecs}{\end{corollary}}
\newcommand{\eers}{\end{exercise}}
\newcommand{\eexs}{\end{example}}
\newcommand{\eems}{\end{example}}
\newcommand{\els}{\end{lemma}}
\newcommand{\eles}{\end{lemmaex}}
\newcommand{\ets}{\end{theorem}}
\newcommand{\eds}{\end{definition}}
\newcommand{\eps}{\end{proposition}}

\newcommand{\bi}{\begin{itemize}}
\newcommand{\ei}{\end{itemize}}
\newcommand{\ben}{\begin{enumerate}}
\newcommand{\een}{\end{enumerate}}

\def\vbar{\mathchoice{\vrule height6.3ptdepth-.5ptwidth.8pt\kern-.8pt}
   {\vrule height6.3ptdepth-.5ptwidth.8pt\kern-.8pt}
   {\vrule height4.1ptdepth-.35ptwidth.6pt\kern-.6pt}
   {\vrule height3.1ptdepth-.25ptwidth.5pt\kern-.5pt}}
\def\fudge{\mathchoice{}{}{\mkern.5mu}{\mkern.8mu}}
\def\bbc#1#2{{\rm \mkern#2mu\vbar\mkern-#2mu#1}}
\def\bbb#1{{\rm I\mkern-3.5mu #1}}
\def\bba#1#2{{\rm #1\mkern-#2mu\fudge #1}}
\def\bb#1{{\count4=`#1 \advance\count4by-64 \ifcase\count4\or\bba A{11.5}\or
   \bbb B\or\bbc C{5}\or\bbb D\or\bbb E\or\bbb F \or\bbc G{5}\or\bbb H\or
   \bbb I\or\bbc J{3}\or\bbb K\or\bbb L \or\bbb M\or\bbb N\or\bbc O{5} \or
   \bbb P\or\bbc Q{5}\or\bbb R\or\bbc S{4.2}\or\bba T{10.5}\or\bbc U{5}\or
   \bba V{12}\or\bba W{16.5}\or\bba X{11}\or\bba Y{11.7}\or\bba Z{7.5}\fi}}

\def \A {{\cal{A}}}

\def \GG {{\cal{G}}}
\def \C {{\cal{C}}}

\def\sqr#1#2{{\vcenter{\vbox{\hrule height .#2pt
                             \hbox{\vrule width .#2pt height#1pt \kern#1pt
                                   \vrule width .#2pt}
                             \hrule height .#2pt}}}}

\def\pmb#1{\setbox0=\hbox{#1}%
   \kern-.025em\copy0\kern-\wd0
   \kern.05em\copy0\kern-\wd0
   \kern-.025em\raise.0433em\box0 }
\def\sqr#1#2{{\vcenter{\vbox{\hrule height.#2pt
     \hbox{\vrule width.#2pt height#1pt \kern#1pt
   \vrule width.#2pt}\hrule height.#2pt}}}}

\def\B{{\mathbb B}}
\def\N{{\mathbb N}}
\def\Z{{\mathbb Z}}
\def\R{{\mathbb R}}

\def\P{{\mathbb P}}

\def\bs{\backslash}

\def\reff#1{(\ref{#1})}
\def \ind {\hbox{ 1\hskip -3pt I}}
\newcommand {\cro}[1] {\left[ {#1} \right]}
\newcommand {\acc}[1] {\left\{ {#1} \right\}}
\newcommand {\pare}[1] {\left( {#1} \right)}

\newcommand {\sous}[1] {\underline{#1}}

\begin{document}
\title{Lower bounds on fluctuations for internal DLA }
\author{
  \renewcommand{\thefootnote}{\arabic{footnote}}
  Amine Asselah\footnotemark[1]
  \and
  \renewcommand{\thefootnote}{\arabic{footnote}}
  Alexandre Gaudilli\`ere\footnotemark[2]}
\date{}

\footnotetext[1]{
    LAMA, Universit\'e Paris-Est -- e-mail: amine.asselah@univ-paris12.fr
  }
\footnotetext[2]{
  LATP, Universit\'e de Provence, CNRS,
  39 rue F. Joliot Curie,
 13013 Marseille, France\\
  \indent\indent e-mail: gaudilli@cmi.univ-mrs.fr
}
\pagestyle{myheadings}
\markboth
    {Fluctuations for internal DLA}
    {Fluctuations for internal DLA}

\maketitle
\begin{abstract}
We consider internal diffusion limited aggregation in dimension
larger than or equal to two. This is a random cluster growth model,
where random walks start at the origin of the $d$-dimensional lattice,
one at a time, and stop moving when reaching a site that is not
occupied by previous walks.

When $n$ random walks are sent from the origin,
we establish a lower bound for the inner and outer errors
fluctuations of order square root of the logarithm of $n$.
When dimension is larger or equal to three, this lower bound
matches the upper bound recently obtained in independent
works of \cite{AG2} and \cite{JLS2}.
Also, we produce as a corollary of our proof of \cite{AG2},
an upper bound for the fluctuation of the inner error in a
specified direction.
\end{abstract}

\section{Introduction}
We establish lower bounds for the inner and
outer errors in internal diffusion limited aggregation (internal DLA).
Internal DLA models a discrete cluster growth, and is defined as follows.
Let $\Lambda$ be a subset of $\Z^d$ which represents 
the explored region at time 0. Let $N$ be an integer,
and $\xi=(\xi_1,\dots, \xi_N)$ be the initial positions 
of $N$ independent simple random walks $S_1,\dots,S_N$ on $\Z^d$.
The cluster of volume $N$, denoted $A(\Lambda, \xi)$ 
is obtained inductively as follows. First, $A(\Lambda,0)=\Lambda$.
Now, assume $A(\Lambda,k-1)$ is obtained, and define
\be{time-settling}
\tau_k=\inf\acc{t\ge 0:\ S_k(t)\not\in A(\Lambda,k-1)},\quad\text{and}\quad
A(\Lambda,k)=A(\Lambda,k-1)\cup \{S_k(\tau_k)\}.
\ee
We call explorers the random walks obeying the aggregation rule
\reff{time-settling}. We say that the explorer $k$ settles
at time $\tau_k$. When $\Lambda=\emptyset$, we also denote
$A(\emptyset,\xi)$ by $A(\xi)$.

In dimensions two and more,
Lawler, Bramson and Griffeath~\cite{lawler92} prove that in order
to cover, without holes, a sphere of radius $n$, we need about the
number of sites of $\Z^d$ in this sphere. In other words,
the asymptotic shape of the cluster
is a sphere. Then, Lawler in~\cite{lawler95} shows a
subdiffusive upper bound for the worse fluctuation to the
spherical shape. More precisely,
the latter result is formulated in terms of inner and outer errors,
which we now introduce with some notation.
We denote with $\|\cdot\|$ the euclidean norm on $\R^d$.
For any $x$ in $R^d$ and $r$ in $\R$, set
\be{ball-dfn}
B(x,r) = \left\{ y\in\R^d :\: \|y-x\| < r \right\}
\quad\mbox{and}\quad \B(x,r) = B(x,r) \cap \Z^d.
\ee
For $\Lambda\subset\Z^d$,
$|\Lambda|$ denotes the number of sites in $\Lambda$.
The inner error $\delta_I(n)$ is such that
\be{def-inner}
n-\delta_I(n)=\sup\acc{r \geq 0:\: \B(0,r)\subset A\big(|\B(0,n)|\big)}.
\ee
Also, the outer error $\delta_O(n)$ is such that
\be{def-outer}
n+\delta_O(n)=\inf\acc{r \geq 0:\: A\big(|\B(0,n)|\big)\subset \B(0,r)}.
\ee
Note that $\delta_I$ and $\delta_O$ represent the worse fluctuations
to the spherical shape.

In 2010, Asselah-Gaudilli\`ere in \cite{AG1,AG2}, 
and Jerison-Levine-Sheffield in \cite{JLS1,JLS2},
independently showed the following upper bound.
\bt{theo-upper}
When dimension $d\ge 3$,
there are constants $\{\beta_d,\ d\ge 3\}$
such that with probability~1,
\be{theo-res2}
\limsup_{n\to\infty} \frac{\delta_I(n)}{\sqrt{\log(n)}}\le \beta_d,
\quad\text{and}\quad
\limsup_{n\to\infty} \frac{\delta_O(n)}{\sqrt{\log(n)}}\le \beta_d.
\ee
When dimension is $2$,
there is a constant $\{\beta_2\}$ such that with probability~1,
\be{theo-res1}
\limsup_{n\to\infty} \frac{\delta_I(n)}{\log(n)}\le \beta_2,
\quad\text{and}\quad
\limsup_{n\to\infty} \frac{\delta_O(n)}{\log(n)}\le \beta_2.
\ee
\et
All these studies are concerned with upper bounds on
$\delta_I$ and $\delta_O$, and the matter of showing that
they were indeed realized remained untouched.
We present now two results on lower bounds. The first one is
independent from previous results.
This bound which characterizes the worse fluctuations
is optimal in $d\ge 3$. Morally, it says that if there are no
{\it deep} holes, then {\it long} tentacles form.
\bp{prop-simple}
When dimension is three or more, set $h(k)= \sqrt{\log(k)}$
there is $\alpha_d>0$ such that
\be{prop-d3}
\lim_{n\to\infty} P\pare{ \exists k>n,\  \delta_O(k)\ge
\alpha_d h(k)\quad\big|\quad
 \delta_I(n)< \alpha_dh(n)}=1.
\ee
When dimension is two, set $h(k)= \sqrt{\log(k)\log(\log(k))}$.
There is $\alpha_2>0$ such that
\be{prop-d2}
\lim_{n\to\infty} P\pare{ \exists k>n,\  \delta_O(k)\ge
\alpha_2 h(k)\quad\big|\quad
 \delta_I(n)< \alpha_2h(n)}=1.
\ee
\ep
\br{rem-io}
Note that as a consequence of \reff{prop-d3}, we have
$\alpha_d$ positive such that
\be{main-io}
P\pare{ \delta_O(n)\ge \alpha_d\sqrt{\log(n)} \text{ or } 
\delta_I(n)\ge \alpha_d\sqrt{\log(n)},\ \text{ i.o.}}=1,
\ee
\er
The second result uses the
upper bound of Theorem~\ref{theo-upper}. It says 
that {\it deep} holes do form.
\bt{theo-lower}
There are positive constants $\{\alpha_d,d\ge 2\}$ such that
when $d\ge 3$, 
\be{main-A}
P\big(\delta_I(n)> \sqrt{ \alpha_d \log(n)}\quad i.o.\big)=1
\ee
\et
\br{rem-conj} We believe that in dimension 2, 
for some $\alpha_2>0$, almost surely both $\delta_I(n)$ and 
$\delta_O(n)$ are larger than $\alpha_2 \log(n)$ infinitely often.
However, the way of realizing this event 
is probably different from what
we describe in $d\ge 3$.
\er
We present now a side result dealing with 
fluctuation in a given direction.
The results on fluctuations for internal DLA \cite{lawler92,lawler95,
AG1,AG2,JLS1,JLS2} all focus on inner and outer errors.
The paper \cite{JLS3} if still of a different nature:
it addresses averaged error, in a sense
where inner and outer errors can cancel each other out.

Though interesting, 
fluctuations in a given direction remained untreated.
This is presented
as an open problem in Section 3 of \cite{JLS2}, and we note that
the proof of \cite{AG2} yields the following
bound for the inner error in a given direction.
\bp{prop-AG}
There are positive constants $\{\kappa_d,d\ge 2\}$ such
that for $z$ be such that $\|z\|<n$, we have in $\Z^d$
\be{AG-directional}
P\pare{ z\not\in A(N_n)}\le
 \left\{ \begin{array}{ll}
\exp\pare{-\kappa_2 \frac{(n-\|z\|)^2}{\log(n)}}& \mbox{ for } d=2 \, , \\
\exp\pare{-\kappa_d (n-\|z\|)^2}& \mbox{ for } d\ge 3 \, .
\end{array} \right.
\ee
\ep
The rest of the paper is organized as follows.
In Section~\ref{sec-notation}, we set notation and recall
useful results. In Section~\ref{sec-outer}, we deal with
the outer error and prove Proposition
\ref{prop-simple}. In Section~\ref{sec-inner}, we deal with
the inner error, and prove Theorem~\ref{theo-lower}. Finally,
in Section~\ref{sec-last}, we explain how to read
Proposition~\ref{prop-AG} as a corollary of our previous work
\cite{AG2}.

\section{Notation and Prerequisite}\label{sec-notation}

Let $S:\N\to\Z^d$ denotes a simple random
walk on $\Z^d$. When the initial condition
is $S(0)=z\in \Z^d$, its law is denoted $\P_z$.
The first time $S$ hits a domain $\Lambda\subset \Z^d$, 
is denoted $H(\Lambda)$, or $H(S;\Lambda)$ to emphasize
that $S$ is the walk. 

For a positive $\gamma$,
we denote by $\rho(\gamma)$ the radius of the largest
ball centered at 0 whose volume is less than $\gamma$. In other words,
\be{not1}
\rho(\gamma)=\sup\acc{n\ge 0: |\B(0,n)|\le \gamma}.
\ee
The abbreviation $b(n)=|\B(0,n)|$ is also handy. 

We now consider
a configuration $\eta\in \N^{\Z^d}$ with a finite number
of particles 
\be{not3}
|\eta|:=\sum_{z\in \Z^d} \eta(z)<\infty.
\ee
For $R>0$, $\Lambda\subset \Z^d$, and
$z\in \B(0,R)$ (resp. $z\in \partial \B(0,R)$), 
we denote by
$W_R(\Lambda;\eta,z)$ the number of explorers visiting
$z$ {\it before} they exit $\B(0,R)$ (resp. {\it at the moment}
they reach $\partial \B(0,R)$), when
the explorers start on $\eta$
with an explored region $\Lambda$.
The positions of explorers settled before they exit
$\B(0,R)$ is denoted $A_R(\Lambda;\eta)$, which is also the
support of $W_R(\Lambda;\eta,z)$ in $\B(0,R)$. 
Note that $A_R(\Lambda;\eta)$ is contained in $\B(0,R)$. 

When we build the
internal DLA cluster from $n$ independent random walks,
we also consider the unrestricted trajectories. This gives
a natural coupling between explorers and independent random walks.
Note that if $\B(0,R)\subset \Lambda$ and $z\in \partial \B(0,R)$, 
then 
\[
W_R(\Lambda;\eta,z)=W_R(\B(0,R);\eta,z),
\]
and this corresponds
to the number of independent random walks which exit
$\B(0,R)$ at site $z$. For simplicity, we call this number
$M_R(\eta,z)$. In general, under a natural coupling, 
for $z\in \B(0,R)\cup\partial \B(0,R)$, 
\be{lawler-3}
W_R(\Lambda, \eta,z) \le M_R(\eta,z).
\ee

\subsection{The abelian property}\label{sec-abelian}
Diaconis and Fulton~\cite{diaconis-fulton}
allow explorers to start on distinct sites, and
show that the law of the cluster is invariant under permutation of
the order in which explorers are launched. This invariance is
named {\it the abelian property}. As a consequence,
one can realize the cluster by sending many {\it exploration waves}.
Let us illustrate this observation by building
$A(\emptyset,(n+m)\delta_0)$ in three waves, since 
we need later this very example. The first wave
consists in launching $n$ explorers,
and they settle in $A_1=A(\emptyset,n\delta_0)$, which we call
the cluster {\it after the first wave}. Then, we launch $m$
explorers that we color {\it green} for simplicity, 
and if they reach $\partial \B(0,R)$ before settling,
then we stop them on $\partial \B(0,R)$.
The settled green explorers make up the cluster $A_R(A_1; m\delta_0)$.
The cluster after the second wave is then
\[
A_2=A_1\cup A_R(A_1; m\delta_0).
\]
For $z\in \partial \B(0,R)$, we call $\zeta_R$
the configuration of the green explorers 
stopped on $\partial \B(0,R)$. In other words,
\[
\zeta_R=\acc{W_R(A_1;m\delta_0,z),\ z\in \partial \B(0,R)}.
\]
Then, the cluster after the third wave is obtain as we launch
the stopped green explorers, and
\[
A_3=A_2\cup A(A_2;\zeta_R).
\]
The abelian properties implies that $A(\emptyset,(n+m)\delta_0)$
equals in law to $A_3$. It is convenient to think of the growing cluster
as evolving in discrete time,
where time counts the number of exploration waves.

\subsection{On the harmonic measure.}\label{sec-harmonic}
We first recall a well known property of Poisson variables.
\bl{lem-obvious}
Let $\{U_n,n\in \N\}$ be an i.\,i.\,d.\, sequence with value 
in a set $E$, and $\{E_1,\dots,E_n\}$ a partition of $E$.
Then, if $X$ is an independent Poisson random variable
of parameter $\lambda$, and if 
\[
X_i=\sum_{n\le X} \ind_{U_n\in E_i},
\]
then $\{X_i,\ i=1,\dots,n\}$ are independent Poisson
variables with $E[X_i]=\lambda\times P(U\in E_i)$.
\el
Now, for $z\in \Z^d\bs\{0\}$, let
\be{not2}
\Sigma(z)=\partial \B(0,\|z\|),
\ee
and note that $z$ belongs to $\Sigma(z)$ since there
is $z'\in \B(0,\|z\|)$ and $\|z-z'\|=1$
(see Lemma 2.1 of \cite{AG1}). 
When we start independent random walks on $\eta$, and
for $z\not= 0$, $h>0$, and $\Lambda\subset \Sigma(z)$, we denote
by $N_z(\eta,\Lambda,h)$ the number of walks initially on $\eta$
which hit $\Sigma(z)$ on $\Lambda$, and then visit $z$ before
exiting $\B(0,\|z\|+h)$.
As a consequence of Lemma~\ref{lem-obvious}, we have
\bc{cor-obvious}
Let $X$ be a Poisson variable of parameter $\lambda$, and
$\eta=X\delta_0$. For any $z\not= 0$, $h,h'>0$, and $\Lambda,\Lambda'
\in \Sigma(z)$ with $\Lambda\cap\Lambda'=\emptyset$, we have
that  $N_z(\eta,\Lambda,h)$ and  $N_z(\eta,\Lambda',h')$ are independent
Poisson variable, and
\be{mean-N}
\begin{split}
E\cro{N_z(\eta,\Lambda,h)}=
&\lambda\times \P_0\Big(S\big(H(\Sigma(z))\big)\in \Lambda,
\ H(z)<H\big(\B^c(0,\|z\|+h)\big)\Big)\\
=&\lambda\sum_{y\in \Lambda}\P_0\Big(S\big(H(\Sigma(z))=y\big)\Big)
\P_y\Big(\ H(z)<H\big(\B^c(0,\|z\|+h)\big)\Big).
\end{split}
\ee
\ec
By combining well known asymptotics of the harmonic measure
with Corollary~\ref{cor-obvious}, we obtain the following lemma.
\bl{lem-main} 
Assume that dimension is two or more.
Let $\eta$ be as in Corollary~\ref{cor-obvious}.
For $z\not= 0$, and $R>0$ let $\Lambda=\B(z,R)\cap \Sigma(z)$, and
$\Lambda'=\Sigma(z)\bs \Lambda$. There is $\kappa>0$, independent
of $z$, and $R$ such that
\be{not4}
P\big(N_z(\eta,\Lambda,\infty)=0,\ N_z(\eta,\Lambda',R)=0\big)\le
\exp\pare{-\kappa \frac{\lambda R}{\|z\|^{d-1}}}.
\ee
\el

\bpr
Since $N_z(\eta,\Lambda,\infty)$ and  $N_z(\eta,\Lambda',R)$
are independent Poisson variables, we have
\be{not5}
P\big(N_z(\eta,\Lambda,\infty)=0,\ N_z(\eta,\Lambda',R)=0\big)=
\exp\big(-E\cro{N_z(\eta,\Lambda,\infty)}-E\cro{
N_z(\eta,\Lambda',R)}\big).
\ee
It remains to compute expected values.
To estimate $E\cro{N_z(\eta,\Lambda,\infty)}$, we recall that
there is a constant $c_d$ such that
\be{not6}
\P_y(H(z)<\infty)\le \frac{c_d}{1+\|y-z\|^{d-2}}.
\ee
Using \reff{mean-N}, there is a constant $c$
\be{not7}
\begin{split}
E\cro{N_z(\eta,\Lambda,\infty)}\le &
\sum_{y\in \Lambda} \lambda\P_0\Big(S\big(H(\Sigma(z))\big)=y\Big)\times
\frac{c_d}{1+\|y-z\|^{d-2}}\\
\le & \frac{c}{2}\frac{\lambda}{\|z\|^{d-1}}
\Big(1+\sum_{k=1}^{R} \frac{k^{d-2}}{1+k^{d-2}}\Big)\le c
\frac{\lambda R}{\|z\|^{d-1}}.
\end{split}
\ee
Now, to estimate $E\cro{N_z(\eta,\Lambda',R)}$, we recall
Lemma 5(b) of \cite{JLS1}, which states that for some
constant $c_d'$, for $y\in \Sigma(z)$
\be{not8}
\P_y\pare{ H(z)< H\big( \B^c(0,\|z\|+R)\big)}\le
\frac{ c_d' R^2}{\|z-y\|^d}
\ee
Using \reff{mean-N} and \reff{not8}, there is a constant $c'$ such that
\be{not9}
\begin{split}
E\cro{N_z(\eta,\Lambda',R)}\le &\sum_{y\in \Sigma(z)\bs \B(z,R)}
\lambda \P_0\Big(S\big(H(\Sigma(z))\big)=y\Big) 
\frac{ c_d' R^2}{\|z-y\|^d}\\
\le &\  \frac{c'}{2} \frac{\lambda R^2}{\|z\|^{d-1}}
\sum_{k=R}^{2\|z\|+2} \frac{k^{d-2}}{k^d}\le \ 
c' \frac{\lambda R}{\|z\|^{d-1}}
\end{split}
\ee
Combining \reff{not7}, and \reff{not9}, we obtain the desired result.
\epr
\section{The outer error}\label{sec-outer}
In this section, we prove Proposition~\ref{prop-simple}.
Let us explain the proof
in dimension three or more, and explain in Remark~\ref{rem-d2}
of Step 2 how we adapt the proof to dimension two.

For positive reals $\alpha$ and $\gamma$, to be chosen later,
we set $h(n)=\alpha\sqrt{\log(n)}$ and $\sous{L}(n)=
\gamma\sqrt{\log(n)}$ for estimates on the outer and inner
fluctuations. 
Even though we eventually take $h(n)=\sous{L}(n)$,
it is useful to keep in mind their distinct nature.
The limit \reff{main-io} follows if for some
small $\gamma=\alpha$, we have
\be{main-1}
\lim_{n\to\infty} P\pare{ \exists k\ge n,\ \delta_O(k)\ge h(k)
\ \big|\ \delta_I(n)<\sous{L}(n)}=1.
\ee
Indeed,
\be{main-2}
\begin{split}
P\big( \exists k\ge n,&\ \delta_O(k)\ge h(k)\quad \text{ or }\quad
\delta_I(k)\ge \sous{L}(k)\big)\\
\ge& P\big(\delta_I(n)\ge \sous{L}(n),\text{ or }
\  \exists k\ge n,\ \delta_O(k)\ge h(k) \big)\\
\ge &P\big(\delta_I(n)\ge \sous{L}(n)\big)
+P\pare{ \exists k\ge n,\ \delta_O(k)\ge h(k)\ \big|\ 
\delta_I(n)< \sous{L}(n)}P\big(\delta_I(n)< \sous{L}(n)\big)\\
\ge &1-\Big(1-P\big( \exists k\ge n,\ \delta_O(k)\ge 
h(k)\ \big|\ \delta_I(n)< \sous{L}(n)\big)\Big)
P\big(\delta_I(n)< \sous{L}(n)\big).
\end{split}
\ee
We now prove \reff{main-1}. For integer $n$, assume that
$A(b(n)\delta_0)$ is realized. 
Let $X_n$ be a Poisson random variable with
parameter $\lambda_n= |\B(0,n+h(n))\bs\B(0,n)|$.
We realize the cluster $A((b(n)+X_n)\delta_0)$ through three
exploration waves, as explained in Section~\ref{sec-abelian}.
After the first wave with $b(n)$ explorers,
we launch $X_n$ explorers, the green ones, and stop them on
$\Sigma:=\partial \B(0,n-\sous{L}(n))$. Under the event
$\{\delta_I(n)<\sous{L}(n)\}$, note that for $z\in \Sigma$
\[
W_{n-\sous{L}(n)}(A(b(n)\delta_0);X_n\delta_0,z)=
M_{n-\sous{L}(n)}(X_n\delta_0,z).
\]
The configuration of stopped random walks on $\Sigma$
is denoted $\zeta$. Note that $\zeta$ is independent of
$A(b(n)\delta_0)$. For $z\in \Sigma$, we call cov$(z)$ the event that
the $\zeta(z)$ green explorers starting on $z$
produce a cluster $A(\emptyset,\zeta(z)\delta_z)$
which satisfies
\be{one-finger}
A(\emptyset,\zeta(z)\delta_z)\cap \B^c(0,n+4h(n))\not= \emptyset.
\ee
Note that the green explorers contributing to cov$(z)$ start
on the positions of the random walks stopped on $z\in\Sigma$,
which are associated with the green explorers.
Assume, for a moment that when cov$(z)$ happens, there is a tentacle of 
$A((b(n)+X_n)\delta_0)$ which protrudes $\B(0,n+4h(n))$.
Assume also that under condition on $X_n$, we have
\be{cond-Xn}
n+4h(n)\ge R_n+h(R_n),\quad\text{where}\quad
R_n=\rho(b(n)+X_n) .
\ee
We would deduce that $\delta_O(R_n)\ge h(R_n)$.
 
We now proceed through four steps. First, we show that \reff{one-finger}
implies that the final cluster is not inside $\B(0,n+4h(n))$.
Secondly, we estimate the cost of producing a tentacle realizing
cov$(z)$. Then, we establish conditions ensuring \reff{cond-Xn}.
Finally, we show that for an appropriate choice of $\alpha$,
one event cov$(z)$ realizes for some $z\in \Sigma$.
\paragraph{Step 1: Coupling.}
By coupling, it is easy to see that for any subset $\Lambda$,
and $z\in \Sigma$
\be{step-3}
A(\emptyset;\zeta(z)\delta_z)\subset
\Lambda\cup A(\Lambda;\zeta(z)\delta_z)\subset
\Lambda\cup A(\Lambda;\zeta).
\ee
If we denote by $A_2$ the cluster after the second
exploration wave (see Section~\ref{sec-abelian}), then
we have with an equality in law
\be{step-4}
A_2\cup A(A_2; \zeta)=
A\big((b(n)+X_n)\delta_0\big).
\ee
Now, using \reff{step-3}, \reff{step-4}, 
and \reff{one-finger}, we conclude that 
the final cluster is not in $\B(0,n+4h(n))$.

\paragraph{Step 2: Long tentacles.}
To produce cov$(z)$, we first bring a number 
of green explorers at $z$ proportional to $h(n)$, and force them to make
a tentacle normal to $\Sigma$ at $z$, with a height $4h(n)+
\sous{L}(n)$. More precisely, draw unit cubes centered on the 
points of the sequence
\[
x_n=\big(\|z\|+n\big) \frac{z}{\|z\|}\in \R^d\quad(\text{and}\quad
\|x_n\|=\|z\|+n).
\]
Each such cube contains at least a site of $\Z^d$, say $z_n$.
Note that $\|z_n-z_{n-1}\|\le 2 \sqrt d+1 $, so that we can
exhibit a sequence $\{z=y_1,y_2,\dots,y_N\}$ of nearest neighbors
in $\Z^d$ such that $\|y_N-z\|\ge 4h(n)+\sous{L}(n)$, with
$N\le c(4h(n)+\sous{L}(n))$ for some constant $c$ independent of $n$.
Now, if $M_{n-\sous{L}(n)}(X_n\delta_0,z)\ge N$, and if we launch
the green explorers stopped on $z$ and force them to walk along
the sequence $\{y_1,y_2,\dots,y_N\}$, with the $k$-th explorer
settling on $y_k$, then we realize cov$(z)$, and its probability
is larger than
\[
\pare{\frac{1}{2d}}^{\sum_{k=1}^N k}\ge 
\exp\big(-c (4h(n)+\sous{L}(n))^2\big).
\]
\br{rem-d2} In dimension 2, there is a better strategy
to build a tentacle. Since we believe that it yields an estimate
which is not optimal, we do not give the full proof, 
but give enough details of the construction so that the interested
reader can easily fill the details. 

We first bring a larger number of green
explorers at $z$, about $ \frac{h(n)}{\log^2(h(n))}\times h(n)$.
The probability of so doing is larger than 
\be{forcing-many}
\exp\Big(-C \frac{h(n)}{\log^2(h(n))}\log\big(\frac{h(n)}{\log^2(h(n))}
\big) \times h(n)\Big)\ge \exp\Big(-c\frac{h^2(n)}{\log(h(n))}\Big),
\ee
for two positive constants $C,c$. Then, the explorers are
forced to fill sequentially cylindrical compartments
of a telescope-like domain that we now describe. Let $R$ be the integer
part of $4h(n)+\sous{L}(n)$, and
divide $B(z,R)$ into $R$ shells of length $h_1,\dots,h_R$ with
for $i=1,\dots,R$
\be{def-hi}
h_i=\frac{R}{i\times \log(R)}.
\ee
Choose the sequence of $R$ points of $\R^d$
\[
\forall i\in \{1,\dots,R\}\quad 
x_i=\big(\|z\|+h_i\big) \frac{z}{\|z\|}
\]
There is $z_i\in \partial \B(z,h_i)$ with $\|z_i-x_i\|\le 2$.
Now, assume we have brought $n_i:=A\pi h_{i+1}^2$ (for a fixed
constant $A$) explorers in $\B(z_i,h_i/4)$. 
Then, with a probability larger than $1-1/h_i^2$,
they cover the ball $\B(z_i,h_{i+1})$ by Lemma 1.3 of
\cite{AG2}. Now, if we bring additional explorers in $\B(z_i,h_i/4)$,
they can reach $\B(z_{i+1},h_{i+1}/4)$ with a positive
probability, say $\exp(-\kappa)$: 
indeed, they only need to escape a cube-like domain
centered on $z_i$, of side-length $h_{i+1}/4$ on the side which
contains $z_{i+1}$. The cost of the scenario, for which we
only described the $i$-th step, is therefore of order
\be{scenario-d2}
\prod_{i=1}^{R-1}\pare{1-\frac{1}{h_i^2}}\times e^{-\kappa
\sum_{i=2}^R n_i}\times e^{-\kappa
\sum_{i=3}^R n_i}\dots\times e^{-\kappa
n_R}\ge C\exp\big(-\kappa \sum_{i=2}^R (i-1)n_i\big),
\ee
for positive constants $\kappa,C$. Note that with the choice
of $h_i$ in \reff{def-hi}, and $n_i=A\pi h_{i+1}^2$, we have
with $\sous{L}(n)=h(n)$ and for a positive constant $\kappa'$
\[
\sum_{i=2}^R (i-1)n_i\ge \kappa' \frac{h^2(n)}{\log(h(n))}.
\]
Thus, the probability of bringing $h^2(n)/\log^2(n)$ explorers
in $z$, and the probability of building a tentacles of height $5h(n)$
is larger than 
\be{order-d2}
\exp\Big(-\kappa_2\frac{h^2(n)}{\log(h(n))}\Big),
\ee
for some positive constant $\kappa_2$.
\er

\paragraph{Step 3: Bounding $X_n$.}
We impose that $X_n\le 2\lambda_n$.
\be{ie35}
\begin{split}
X_n &\le 2\big|\B(0,n+h(n))\bs \B(0,n)\big|\\
&\le \big|\B(0,n+2h(n))\bs \B(0,n)\big|
\Longrightarrow X_n+b(n)\le \big|\B(0,n+2h(n))\big|
\Longrightarrow R_n\le n+2h(n).
\end{split}
\ee
The conclusion of \reff{ie35} implies also that $h(R_n)\le 2h(n)$,
and this implies \reff{cond-Xn}.
Note that since $X_n$ is a Poisson variable of mean $\lambda_n$,
there is a constant $c$, such that
\be{excess-1}
P\big(X_n> 2\lambda_n\big)\le \exp\big(-ch(n)n^{d-1}\big).
\ee
\paragraph{Step 4: Many possible tentacles.}
We summarize Step 1 to Step 3, as establishing that
\be{ie26}
\begin{split}
\acc{ X_n\le 2 \lambda_n,\ \exists z\in \Sigma\quad \text{cov}(z),\ 
\delta_I(n)<\sous{L}(n)}\subset&
\acc{ \delta_O(R_n)\ge h(R_n),\ \delta_I(n)<\sous{L}(n)}\\
\subset&\acc{\exists k>n, 
\delta_O(k)\ge h(k),\ \delta_I(n)<\sous{L}(n)}.
\end{split}
\ee
The key observation now is that $\{ \text{cov}(z),\ z\in \Sigma\}$ 
are independent Poisson variables which are also independent
of $\delta_I(n)$. Indeed, $\{ \text{cov}(z),\ z\in \Sigma\}$
deals with the walks associated with the green explorers,
whereas $\delta_I(n)$ depends on the $b(n)$ explorers which
we launch first. Taking probability
on both sides of \reff{ie26}, and dividing by $P(\delta_I(n)<\sous{L}(n))$,
we obtain
\be{ie27}
P\big(\exists k>n, 
\delta_O(k)\ge h(k)\big| \delta_I(n)<\sous{L}(n)\big)\ge
P\big( X_n\le 2 \lambda_n,\ \exists z\in \Sigma\quad \text{cov}(z)\big).
\ee
Thus,
\be{ie28}
P\big(\exists k>n,
\delta_O(k)\ge h(k)\big| \delta_I(n)<\sous{L}(n)\big)\ge
P\big(\ \exists z\in \Sigma \text{cov}(z)\big)-
P\big( X_n> 2 \lambda_n\big).
\ee
Now, note that for some positive constants $\kappa,\kappa'$
\be{main-3}
\begin{split}
P\pare{\bigcup_{z\in \Sigma} \text{cov}(z)}=& 1-
\prod_{z\in \partial \B}
\Big(1-P\big(\text{cov}(z)
\cap\ \{\zeta(z)>c\big(4h(n)+\sous{L}(n)\big)\}\big)\Big)\\
\ge& 1-\exp\pare{-\kappa' n^{d-1}\exp\pare{-\kappa(h^2(n)+\sous{L}^2(n))}}.
\end{split}
\ee
Now, there is $\alpha>0$ such that for $h(n)\le \alpha\sqrt{\log(n)}$,
and $\sous{L}(n)=h(n)$
\be{step-7}
\lim_{n\to\infty}n^{d-1}\exp\pare{-\kappa(h^2(n)+h(n))}=\infty.
\ee
We now combine \reff{ie28}, \reff{ie26}, \reff{main-3} and
\reff{step-7} to conclude the proof.
In dimension 2, the estimate \reff{main-3} has to be replaced
with \reff{order-d2}.

\section{The inner error}\label{sec-inner}
In this section, we prove Theorem~\ref{theo-lower}. 
We show that in the process of going
from a cluster of volume $b(n)$ to one of volume $b(2n)$,
chances tend to one as $n$ tends to infinity,
that there appears a cluster $A$
whose inner error is larger than $\alpha \sqrt{\log(\rho(A))}$
for some positive $\alpha$ independent of $n$.
To do so, we launch many exploration waves, each one is
 made up of a Poisson number of explorers.

We proceed inductively.
For positive reals $\alpha,\beta$, to be chosen later,
we set $h(n)=\alpha\sqrt{\log(n)}$, and $\bar{L}(n)
=\beta\sqrt{\log(n)}$. $\bar{L}(n)$ will
refer to an outer radius. 
Let $\{\GG_n,\ n\ge 0\}$ denotes the natural filtration associated with
the evolution by waves. 

First, we launch $b(n)$ explorers.
Assume that explorers of wave $k-1$ have been launched, and are settled. 
Knowing $\GG_{k-1}$, the size of the $k$-th
wave, denoted $X_k$, is a Poisson variable of parameter
\be{ie1}
\lambda(k)=\big|\B(0,R_{k-1}+2h(R_{k-1}))\bs\B(0,R_{k-1}))\big|,
\quad\text{where}\quad R_{k-1}=\rho(b(n)+X_1+\dots+X_{k-1}).
\ee
Since $R_{k-1}$ is of order $n$,
each wave fills approximately a peel of width $2h(n)$, and
$n/2h(n)$ waves fill approximately $\B(0,2n)$. 
We prove in this section that for an appropriate $\alpha$
\be{ie3}
\lim_{n\to\infty} P\pare{\bigcup_{1\le k<n/2h(n)} 
\acc{\delta_I(R_k)> \alpha \sqrt{\log(R_k)}}}=1.
\ee
We now proceed in estimating the probability of observing
a {\it deep} hole after each exploration waves.
We set $\A_k=\acc{\delta_I(R_k)> \alpha \sqrt{\log(R_k)}}$.

\paragraph{On the holes left after wave $k-1$.}
Observe that by definition, on $\A^c_{k-1}$
\[
\B(0,R_{k-1}-h(R_{k-1}))\subset A\big(b(n)+X_1+\dots+X_{k-1}\big),
\]
which implies that
\be{ie4}
\Big(\B(0,R_{k-1})\bs \B(0,R_{k-1}-h(R_{k-1}))\cup
\partial \B(0,R_{k-1})\Big)\cap A\big(b(n)+X_1+\dots+X_{k-1}\big)\not=
\emptyset.
\ee
Choose any $Z_k$ in the intersection of the non-empty set of \reff{ie4},
and note that
\be{ie24}
R_{k-1}-h(R_{k-1})\le \|Z_k\|\le R_{k-1}+1.
\ee
Recall that we have defined $\Sigma(Z_k):=\partial \B(0,\|Z_k\|)$.
We launch the $X_k$ explorers, that we name the {\it green explorers},
and we stop them as they reach $\Sigma(Z_k)$.
The green explorers which settle before reaching $\Sigma(Z_k)$ 
play no role here, and we bound the number of
green explorers stopped on some region $\Lambda\in \Sigma(Z_k)$,
by the number of corresponding random walks exiting
$\Sigma(Z_k)$ on $\Lambda$. Thus, if we choose
\be{ie5}
\Lambda_k=\B(Z_k,\bar L(R_k))\cap \Sigma(Z_k),
\quad\text{and }\quad \Lambda'_k=\Sigma(Z_k)\bs \Lambda_k,
\ee
and if we denote
\be{ie6}
I_k=\acc{N_{Z_k}(X_k\delta_0,\Lambda_k,\infty)=0,\ 
N_{Z_k}(X_k\delta_0,\Lambda'_k,7\bar L(R_{k-1}))=0}
\ee
then, on the event $I_k$, green explorers either 
exit a ball of radius $R_{k-1}-h(R_{k-1})+7\bar L(R_{k-1})$, or
do not visit $Z_k$. In other words,
\[
I_k\subset \acc{R_k-\delta_I(R_k)< \|Z_k\|}\cup
\acc{R_k+\delta_O(R_k)\ge \|Z_k\|+7\bar L(R_{k-1})}.
\]
In order to conclude that $\{\delta_I(R_k)\ge h(R_k)\}$ or
$\{\delta_O(R_k)\ge \bar L(R_{k})\}$, we need to find 
conditions on $X_k$ that guarantee that
\be{cond-Xk}
R_k-\|Z_k\|\ge h(R_k),\quad\text{and}\quad
\|Z_k\|-R_k+7\bar L(R_{k-1})\ge \bar L(R_{k}).
\ee

\paragraph{ Conditions on $X_k$ fulfilling \reff{cond-Xk}.}
We call
\[
\C_k=\acc{2\lambda(k)\ge X_k}\cap \acc{X_k\ge \frac{2}{3}\lambda(k)}.
\]
On the one hand, if $X_k\ge \frac{2}{3}\lambda(k)$, then
\be{ie21}
\begin{split}
X_k &\ge \frac{2}{3}\big|\B(0,R_{k-1}+2h(R_{k-1}))\bs \B(0,R_{k-1})\big|\\
&\ge \big|\B(0,R_{k-1}+\frac{4}{3}h(R_{k-1}))\bs \B(0,R_{k-1})\big|
\Longrightarrow R_k\ge R_{k-1}+\frac{4}{3} h(R_{k-1}).
\end{split}
\ee
On the other hand, if $X_k\le 2 \lambda_k$, then
\be{ie22}
\begin{split}
X_k &\le 2\big|\B(0,R_{k-1}+2h(R_{k-1}))\bs \B(0,R_{k-1})\big|\\
&\le \big|\B(0,R_{k-1}+4h(R_{k-1}))\bs \B(0,R_{k-1})\big|
\Longrightarrow R_k\le R_{k-1}+4h(R_{k-1}).
\end{split}
\ee
Now, for $x$ large enough, the following implication is obvious
\be{obvious-1}
x\le y+4h(y)\Longrightarrow h(x)\le \frac{4}{3}h(y)-1
,\quad\text{and}\quad \bar L(x)\le 2 \bar L(y).
\ee
If $n$ is large enough, \reff{obvious-1} and \reff{ie22} imply that
$h(R_k)\le \frac{4}{3}h(R_{k-1})-1$, which in turn, with \reff{ie21}, 
yields
\be{ie23}
R_k\ge R_{k-1}-1+h(R_k).
\ee
Also, $\bar L(R_k)\le 2 \bar L(R_{k-1})$, $\bar L(R_{k-1})\ge
h(R_{k-1})$, in combination with \reff{ie22} imply that
\be{ie25}
R_{k-1}-h(R_{k-1})+7\bar L(R_{k-1})-\bar L(R_{k})\ge
R_{k-1}+4h(R_{k-1})\ge R_k.
\ee
Thus, if $\C_k\cap \A^c_{k-1}$ holds, then 
\reff{ie24}, \reff{ie23}  and \reff{ie25} imply that
conditions \reff{cond-Xk} holds.

\paragraph{On a deep hole in one shell.}
We choose an integer $k<n/2h(n)$.
We have seen that knowing $\A^c_{k-1}$
\be{ie10}
I_k\cap \C_k\subset \A_k \cup \acc{\delta_0(R_k)\ge \bar L(R_k)}.
\ee
Taking conditional probabilities on both sides of \reff{ie10},
we obtain, 
\be{ie20}
\begin{split}
\ind_{\A^c_{k-1}}
P\big(\A^c_k\cap\C_k\ \big|\GG_{k-1}\big)=&
\ind_{\A^c_{k-1}}\big(P(\C_k)-P(\A_k\cap \C_k\big|\GG_{k-1})\big)\\
\le &
\ind_{\A^c_{k-1}}\Big(
P(\C_k)-P(I_k\cap \C_k\big|\GG_{k-1})+
P\big(\acc{\delta_0(R_k)\ge \bar L(R_k)}|\GG_{k-1}\big)\Big)\\
\le &P(\acc{\delta_0(R_k)\ge \bar L(R_k)}|\GG_{k-1})
+\ind_{\A^c_{k-1}}\big(1-P(I_k\big|\GG_{k-1})\big).
\end{split}
\ee
Now, we invoke Lemma~\ref{lem-main}
with $\|z\|,R$ and $\lambda$ respectively of order
$n,\bar L(n)$, and $h(n) n^{d-1}$. As a consequence, 
we have on $\C_k\cap \A^c_{k-1}$ for a constant $\kappa$, 
\be{ie-main}
\inf_{k\le n/2h(n)} P(I_k\ |\ \GG_{k-1})\ge \exp(-\kappa h(n) \bar L(n)).
\ee
If we denote $N$ for the integer part of $n/2h(n)$, and proceed
inductively, we obtain
\[
\begin{split}
P\big(\cup_{k\le N} \A_k\big)& 
-P\big(\cap_{k\le N} \C_k\big)\ge-
P\big(\forall k\le N,\  \A^c_k\cap\C_k\big)\\
&\ge -E\cro{\ind_{\forall k<N,\ \A^c_k\cap\C_k}
P\big(\A^c_N\cap\C_N\ \big|\GG_{N-1}\big)} \\
&\ge -E\cro{\ind_{\forall k<N,\ \A^c_k\cap\C_k}
\big(P(\acc{\delta_0(R_N)\ge \bar L(R_N)}|\GG_{N-1})
+1-P(I_N\big|\GG_{N-1})}\big)\\
&\ge -P(\acc{\delta_0(R_N)\ge \bar L(R_N)})-(1-
\exp(-\kappa h(n) \bar L(n)))P\pare{\forall k<N,\ \A^c_k\cap\C_k}\\
&\ge -\sum_{k\le N} P(\acc{\delta_0(R_k)\ge \bar L(R_k)})-
(1-\exp(-\kappa h(n) \bar L(n)))^{N}.
\end{split}
\]
Thus,
\be{ie11}
P\big(\cup_{k\le N} \A_k\big)\ge 1-
\sum_{k\le N} \Big(P\big(\acc{\delta_0(R_k)\ge \bar L(R_k)}\big)+
P(\C_k^c)\Big)-(1-\exp(-\kappa h(n) \bar L(n)))^{N}.
\ee
Now, we have established in \cite{AG2}, that for $\beta$ large
enough, the probability
of $\{\delta_0(R_k)\ge \bar L(R_k)\}$ decays faster than any power
in $n$, whereas the fact that $X_k$ is Poisson implies that for some
constant $c$, we have $P(\C_k^c)\le \exp(-c h(n)n^{d-1})$. The last term
on the last display of \reff{ie11} tends to 0 if
\be{cond-h}
\lim_{n\to\infty} \frac{n}{2h(n)}\exp(-\kappa h(n) \bar L(n))=
\infty.
\ee
In dimension 3 or more, \reff{cond-h} holds
for $\alpha$ small enough.

\section{Proof of Proposition~\ref{prop-AG}}\label{sec-last}
The proof is a direct corollary of formula (3.11) of \cite{AG2}.
We consider actually {\it tiles} of size 1, that is site of $\Z^d$.
Inequality (3.8) of \cite{AG2} shows that for some constant $c_d$
(depending only on dimension) and $R=\|z\|<n$, we have
\be{ag2-1}
E[W_R(\emptyset, N\delta_0,z)]\ge c_d (n-\|z\|).
\ee
Inequality (3.10) of \cite{AG2} is written a little differently as
\be{outer-1}
P\pare{ W_{R}(\emptyset, N\delta_0,z)=0}\le \left\{ \begin{array}{ll}
\exp\pare{-\lambda \kappa_2 (n-\|z\|) +
\lambda^2 c'_2 \log(n)}& \mbox{ for } d=2 \, , \\
\exp\pare{-\lambda \kappa_d (n-\|z\|) +
\lambda^2 c'_d }& \mbox{ for } d\ge 3 \, .
\end{array} \right.
\ee
As we optimize \reff{outer-1} in $\lambda>0$, we obtain
\reff{AG-directional}.


\begin{thebibliography}{99}

\bibitem{AG1} Asselah A., Gaudilli\`ere A., preprint 2010, arXiv:1009.2838
{\em From logarithmic to subdiffusive polynomial fluctuations
for internal DLA and related growth models.}

\bibitem{AG2} Asselah A., Gaudilli\`ere A., preprint 2010,
arXiv:1011.4592, {\em Sub-logarithmic fluctuations for internal DLA.}


\bibitem{diaconis-fulton} Diaconis, P.; Fulton, W.
{\em A growth model, a game, an algebra, Lagrange inversion, and
characteristic classes.}
Rend. Sem. Mat. Univ. Politec. Torino 49 (1991), no. 1, 95--119 (1993).


\bibitem{JLS1} Jerison, D.; Levine, L.; Sheffield, S., preprint 2010,
arXiv:1010.2483 . {\em Logarithmic fluctuations for internal DLA.}

\bibitem{JLS2} Jerison, D.; Levine, L.; Sheffield, S., preprint 2010,
arXiv:1012.3453, {\em Internal DLA in Higher Dimensions }

\bibitem{JLS3} Jerison, D.; Levine, L.; Sheffield, S., preprint 2011,
arXiv:1101.0596, {\em Internal DLA and the Gaussian free field}

\bibitem{lawler92} Lawler, G.; Bramson, M.; Griffeath, D.
{\em Internal diffusion limited aggregation.}
Ann. Probab.  20  (1992),  no. 4, 2117--2140.


\bibitem{lawler95} Lawler, G.
{\em Subdiffusive fluctuations for internal diffusion limited aggregation.}
Ann. Probab.  23  (1995),  no. 1, 71--86.

\end{thebibliography}
\end{document}